\documentclass[reqno,11pt]{amsart}
\usepackage{amssymb}
\usepackage{amsmath}
\usepackage{times}
\usepackage{graphicx}
\usepackage{hyperref}

\newtheorem{theorem}{Theorem}[section]

\begin{document}
\title{On the largest Lyapunov exponent for products of Gaussian matrices}
\thanks{e-mail: vladislav.kargin@gmail.com}
\date{June 2014}
\author{Vladislav Kargin}
\maketitle

\begin{center}
\textbf{Abstract}
\end{center}

\begin{quotation}
The paper provides a new integral formula for the largest Lyapunov exponent
of Gaussian matrices, which is valid in the real, complex and
quaternion-valued cases. This formula is applied to derive asymptotic
expressions for the largest Lyapunov exponent when the size of the matrix is
large and compare the Lyapunov exponents in models with a spike and no
spikes.
\end{quotation}

\bigskip

\section{Introduction}

Let $A_{k}$ be random $d$-by-$d$ matrices and let $P_{n}:=A_{n}A_{n-1}\ldots
A_{1}.$ In a seminal paper \cite{furstenberg_kesten60}, Furstenberg and
Kesten showed that $\lim_{n\rightarrow \infty }n^{-1}\log \left\Vert
P_{n}\right\Vert $ exists under some mild conditions on random matrices $%
A_{k}.$ This limit is called \emph{the largest Lyapunov exponent}. In fact,
the limit of matrices $\left( P_{n}^{\ast }P_{n}\right) ^{1/2n}$ is well
defined (under some conditions on $A_{k}$), and the limit has $d$
non-negative real eigenvalues (Oseledec \cite{oseledec68}, Raghunathan \cite%
{raghunathan79}). The logarithms of these eigenvalues are called the \emph{%
Lyapunov exponents} of matrices $A_{i}$. We denote them $\mu _{1}\geq \ldots
\geq \mu _{d}.$

As Lyapunov exponents are related to the behavior of complex dynamical
systems, much research went into the study of their properties, such as
their analyticity with respect to small matrix perturbations and to changes
in the matrix random process (\cite{ruelle79}, \cite{lepage89}, \cite%
{peres92}), the multiplicity of the largest Lyapunov exponent (\cite%
{guivarch_raugi85}), and others (see papers in \cite{cohen_kesten_newman86}%
). Lyapunov exponents have been generalized to infinite-dimensional
operators (\cite{ruelle82}, \cite{kargin08}), and found important
applications in the study of hydrodynamic stability (\cite{ruelle84}) and
random Schroedinger operators (\cite{bougerol_lacroix85}).

In \cite{furstenberg63}, Furstenberg derived an integral formula for the
largest Lyapunov exponent, which in the case of Gaussian matrices reduces to
a multidimensional integral over a known density. However, there is only a
limited number of cases in which this integral has been computed in a closed
form (\cite{cohen_newman84}, \cite{mtw08}), although in some cases, it was
shown that the exponents can be computed efficiently even if one does non
have an explicit formula (\cite{pollicott10a}).

Recently, in \cite{forrester13} Forrester found a formula for all Lyapunov
exponents of complex Gaussian matrices with a general covariance matrix $%
\Sigma $. His method is based on the Harish-Chandra-Itzykson-Zuber
integration formula and cannot be applied to the case when matrices are
real-valued.

In this paper we derive a new formula for the largest Lyapunov exponent of
Gaussian matrices with general $\Sigma $, which is applicable in the real,
complex and quaternion cases.

Recall that a $d$-by-$d$ Gaussian matrix $A$ with covariance matrix $\Sigma $
has the probability density 
\begin{equation*}
P\left( A\right) =c_{\beta }\left( \det \Sigma \right) ^{-d/2}\exp \left[ -%
\frac{\beta }{2}\mathrm{Tr}\left( A^{\ast }\Sigma ^{-1}A\right) \right] ,
\end{equation*}%
where $\beta =1,$ $2,$ or $4$ for real, complex or quaternion matrices, and $%
c_{\beta }$ is a normalization constant. Equivalently, $A$ can be defined as 
$\Sigma ^{1/2}G,$ where $G$ has independent real, complex, or quaternion
entries whose components are real Gaussian variables with variance $1/\beta $%
.

\begin{theorem}
\label{theorem_main_formula} Let $A_{i}$ be independent $d$-by-$d$ Gaussian
matrices with covariance matrix\emph{\ }$\Sigma ,$ and let their entries be
real, complex or quaternion, according to whether $\beta =1,2,$ or $4.$
Assume that the eigenvalues of $\Sigma $ are $\sigma _{i}^{2}=1/y_{i}.$
Then, the following formula holds for the largest Lyapunov exponent of $%
A_{i},$ 
\begin{equation}
2\mu _{1}=-\gamma +\log \left( \frac{2}{\beta }\right) +\int_{0}^{\infty }%
\left[ \boldsymbol{1}_{[0,1]}\left( x\right) -\prod_{i=1}^{d}\left( 1+\frac{x%
}{y_{i}}\right) ^{-\beta /2}\right] \frac{dx}{x},  \label{formula_main3}
\end{equation}%
where $\gamma \approx 0.5772$ is the Euler constant.
\end{theorem}

As an application of this formula, we derive asymptotic expressions for the
largest Lyapunov exponents of high-dimensional matrices. First, we compute
the asymptotic behavior in the case when $d$ is large and $\Sigma $ does not
have spikes so that all eigenvalues of $\Sigma $ are of the same order.
Next, we consider the case when the covariance matrix $\Sigma $ does have a
\textquotedblleft spike\textquotedblright . We consider the particular case
when all except one eigenvalue of $\Sigma $ equal $1/d$, and the exceptional
eigenvalue does not depend on the dimension. We study the asymptotic
behavior of the largest Lyapunov exponent when the dimension $d$ is large,
and find a significant difference relative to the case of no spikes.

More formally, we have the following result for the no-spike case.

\begin{theorem}
\label{theorem_no_spikes}Suppose that $A_{i}^{\left( d\right) }$ are
independent $d$-by-$d$ Gaussian matrices with covariance matrix $\Sigma
_{d}, $ and let the eigenvalues of $\Sigma _{d}$ be $\theta _{i}^{\left(
d\right) }/d$ where $i=1,\ldots ,d.$ Assume that $\theta _{i}^{\left(
d\right) }$ are bounded, $1\leq \theta _{i}^{\left( d\right) }\leq L,$ and
that 
\begin{equation*}
\lim_{d\rightarrow \infty }\mathrm{Tr}\left( \Sigma _{d}\right) \equiv
\lim_{d\rightarrow \infty }\frac{1}{d}\sum_{i=1}^{d}\theta _{i}^{\left(
d\right) }=\lambda .
\end{equation*}%
Then, for the largest Lyapunov exponent of $A_{i}^{\left( d\right) },$ we
have 
\begin{equation*}
\lim_{d\rightarrow \infty }\mu _{1}^{\left( d\right) }=\frac{1}{2}\log
\lambda .
\end{equation*}
\end{theorem}

The result of Theorem \ref{theorem_no_spikes} is in agreement with the free
probability prediction (see, for example, \cite{kargin07}).

\begin{figure}[tbph]
\includegraphics[width=12cm]{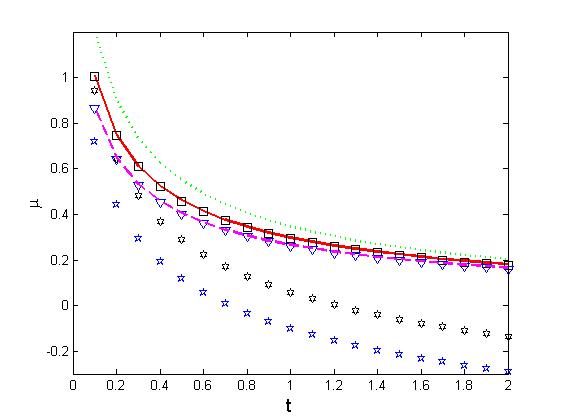} \label{fig:spike}
\caption{The eigenvalues of the covariance matrix equal $1/d$ except one
which equals $1/t$. The blue pentagrams show the largest Lyapunov exponents
when the dimension $d=2$ and $\protect\beta =1$; the black hexagrams are for 
$d=2$ and $\protect\beta =2$; the triangles are for $d=100$ and $\protect%
\beta =1$; the squares are for $d=100$ and $\protect\beta =2$; the red solid
line is the asymptotic prediction for $\protect\beta =2$ and the magenta
dashed line is the prediction for $\protect\beta =1$; the green dotted line
is the free probability prediction $\log (1+1/t)/2$.}
\end{figure}

For the model with a spike we have the following result.

\begin{theorem}
\label{theorem_spikes}Suppose that $A_{i}$ are independent $d$-by-$d$
Gaussian matrices with the covariance matrix $\Sigma _{d},$ and that the
eigenvalues of $\Sigma _{d}$ are $\lambda _{i}^{\left( d\right) }=1/d$ for $%
i=1,\ldots ,d-1,$ and $\lambda _{d}^{\left( d\right) }=1/t$, where $0<t<d.$
In the complex case ($\beta =2$), we have the following estimate$,$ 
\begin{eqnarray}
2\mu _{1} &=&e^{t}\int_{1}^{\infty }e^{-tx}\frac{dx}{x}+O_{t}\left(
1/d\right) ,  \label{formula_asymptotics} \\
&=&-e^{t}\mathrm{Ei}\left( -t\right) +O_{t}\left( 1/d\right) ,  \notag
\end{eqnarray}%
where $\mathrm{Ei}\left( x\right) $ is the exponential integral function. In
the real case ($\beta =1$),%
\begin{equation*}
2\mu _{1}=e^{t/2}\int_{1}^{\infty }e^{-tx/2}\frac{dx}{\sqrt{x}\left( \sqrt{x}%
+1\right) }+O_{t}\left( 1/d\right) ,
\end{equation*}
\end{theorem}

A comparison of the results for the spike and no-spike models shows that a
spike can significantly influence the largest Lyapunov exponent. Indeed, in
the situation of Theorem \ref{theorem_spikes} with $\beta =2$, the formula
in Theorem \ref{theorem_no_spikes} would incorrectly predict that for large $%
d$ the first Lyapunov exponent is close to $\frac{1}{2}\log (1+1/t)$
(because $\lim_{d\rightarrow \infty }\frac{1}{d}\sum_{i=1}^{d}\theta
_{i}^{\left( d\right) }=1+\frac{1}{t}$) instead of the correct $-\frac{1}{2}%
e^{t}\mathrm{Ei}\left( -t\right) $.

Theorem \ref{theorem_spikes} concerns the case when the spike (i.e., the
ratio of the exceptional eigenvalue to the bulk eigenvalues) grows at the
same rate as the dimension. For other regimes, we have the following results.

\begin{theorem}
\label{theorem_no_spikes_2nd}Suppose that $A_{i}$ are independent $d$-by-$d$
complex Gaussian matrices with the covariance $\Sigma ,$ and that the
eigenvalues of $\Sigma _{d}$ are $\lambda _{i}^{\left( d\right) }=1/d$ for $%
i=1,\ldots ,d-1,$ and $\lambda _{d}^{\left( d\right) }=\theta /d$, where $%
\theta >1.$ Then, for the largest Lyapunov exponent $\mu _{1},$ we have 
\begin{equation*}
2\mu _{1}=\frac{\theta -\frac{3}{2}}{d}+O\left( \left( \frac{\theta -1}{d}%
\right) ^{2}\right) +O\left( \frac{1}{d^{2}}\right) .
\end{equation*}
\end{theorem}

Note that for $d\rightarrow \infty $ and $\theta $ fixed, this is
essentially the no-spike model with $\lambda =\lim_{d\rightarrow \infty }%
\frac{1}{d}\sum_{i=1}^{d}\theta _{i}^{\left( d\right) }=1.$ Thus, Theorem %
\ref{theorem_no_spikes_2nd} gives a second order correction to the result in
Theorem \ref{theorem_no_spikes}.

Finally, consider the case of an exceptionally large eigenvalue which grows
faster than the dimension $d.$

\begin{theorem}
\label{theorem_spikes_2nd}Suppose that $A_{i}$ are independent $d$-by-$d$
complex Gaussian matrices with the covariance $\Sigma ,$ and that the
eigenvalues of $\Sigma _{d}$ are $\lambda _{i}^{\left( d\right) }=1$ for $%
i=1,\ldots ,d-1,$ and $\lambda _{d}^{\left( d\right) }=\theta >d.$ Then, for
the largest Lyapunov exponent $\mu _{1},$ we have 
\begin{equation*}
2\mu _{1}=\log \theta -\gamma +O(\frac{d}{\theta }\log \theta ).
\end{equation*}
\end{theorem}

In particular, if $\theta $ grows faster than $d^{\alpha }$ with $\alpha >1,$
then the last term is negligible for large $\theta .$ This happens, for
example, if $d$ is fixed and $\theta $ grows. It is remarkable that the
Lyapunov exponent does not depend on the dimension $d$ if the spike in
covariance $\theta $ is much larger than $d$.

The rest of the paper is organized as follows. Section \ref%
{section_background} collects necessary background information about
Lyapunov exponents. Section \ref{section_integral_formula} is devoted to the
proof of the main formula. Section \ref{section_spike} derives asymptotic
expressions for the largest Lyapunov exponent when the size of the matrices
is large. And Section \ref{section_conclusion} is the conclusion.

\section{Background information about Lyapunov exponents}

\label{section_background}

One of the key facts about Lyapunov exponents is that they satisfy the
following relation: 
\begin{equation}
\mu _{1}+\ldots +\mu _{k}=\lim_{n\rightarrow \infty }\frac{1}{n}\log \mathrm{%
Vol}_{k}\left( y_{1}\left( n\right) ,\ldots ,y_{k}\left( n\right) \right)
\label{characterization_Lyapunov}
\end{equation}%
where $y_{i}\left( n\right) =P_{n}y_{i}\left( 0\right) $ and $\left\{
y_{i}\left( 0\right) \right\} _{i=1}^{k}$ is an arbitrary orthonormal vector
system.

For Gaussian matrices this formula can be significantly simplified. Namely,
let $A_{i}$ be independent Gaussian matrices with covariance matrix\textit{\ 
}$\Sigma .$ The crucial observation is that the distribution of $A_{i}^{\ast
}A_{i}$ is invariant relative to the transformation 
\begin{equation*}
A_{i}^{\ast }A_{i}\rightarrow Q^{\ast }A_{i}^{\ast }A_{i}Q,
\end{equation*}%
where $Q$ is an arbitrary orthogonal matrix. This implies that the changes
in the volume of a $k$-dimensional element are independent from step to step
and that their distribution is the same as if they were applied to the
element spanned by the standard basis vectors $\mathbf{e}_{i}.$ Hence, by
the law of large numbers, 
\begin{eqnarray}
\mu _{1}+\ldots +\mu _{k} &=&\mathbb{E}\log \mathrm{Vol}_{k}\left( A_{1}%
\mathbf{e}_{1},\ldots ,A_{1}\mathbf{e}_{k}\right)  \notag \\
&=&\frac{1}{2}\mathbb{E}\log \det \left( G_{k}^{\ast }\Sigma G_{k}\right) =%
\frac{1}{2}\left. \frac{d}{d\mu }\mathbb{E}\left[ \det \left( G_{k}^{\ast
}\Sigma G_{k}\right) ^{\mu }\right] \right\vert _{\mu =0},
\label{basic_formula0}
\end{eqnarray}%
where $G_{k}$ denotes a random $d$-by-$k$ matrix with the identically
distributed standard Gaussian entries. (For details of the argument see \cite%
{newman86a} and \cite{newman86b}.)

While formula (\ref{basic_formula0}) allows one to compute all Lyapunov
exponents, it is essentially a multidimensional integral which can be
computationally demanding. For this reason, it is of interest to obtain a
more explicit way to calculate the Lyapunov exponents.

For real Gaussian matrices and the simplest situation when $\Sigma =\sigma
^{2}I$ and $I$ is the identity matrix, Newman showed in \cite{newman86b} that%
\begin{equation}
\mu _{i}=\frac{1}{2}\left[ \log \left( 2\sigma ^{2}\right) +\Psi \left( 
\frac{d-i+1}{2}\right) \right] ,  \label{Newman}
\end{equation}%
where $\Psi \left( x\right) $ is the digamma function, $\Psi \left( x\right)
:=\left( \log \Gamma \left( x\right) \right) ^{\prime }.$ (At the positive
integer points, $\Psi \left( n\right) =\sum_{k=1}^{n-1}\frac{1}{k}-\gamma ,$
where $\gamma =0.5772\ldots $ is the Euler constant. At half-integers, $\Psi
\left( n+1/2\right) =\sum_{k=1}^{n}\frac{1}{k-1/2}-2\log 2-\gamma .$ The
asymptotic behavior of the digamma function is given by the formula $\Psi
\left( z\right) =\log z-\frac{1}{2z}-\frac{1}{12z^{2}}\left( 1+O\left( \frac{%
1}{z^{2}}\right) \right) $.)

In particular if we normalize $\sigma ^{2}=1/d,$ then for $d=1$ the largest
Lyapunov exponent $\mu _{1}=[-\log 2-\gamma ]/2$ and for $d\rightarrow
\infty ,$ $\mu _{1}=-\frac{1}{2d}+O\left( \frac{1}{d^{2}}\right) .$

In \cite{forrester13}, Forrester established an analogue of formulas (\ref%
{Newman}) for the complex-valued Gaussian matrices. Namely, Forrester showed
that in the case of complex-valued Gaussian matrices with $\Sigma =\sigma
^{2}I$, 
\begin{equation}
2\mu _{i}=\log \sigma ^{2}+\Psi \left( d-i+1\right)
\end{equation}%
(see Proposition 1 in \cite{forrester13} and note that the absence of $1/2$
before $\Psi $ is a typo in the statement of this proposition.) If $\sigma
^{2}=1/d,$ then for $d=1$ the largest Lyapunov exponent $\mu _{1}=-\gamma /2$
and for $d\rightarrow \infty ,$ $\mu _{1}=-\frac{1}{4d}+O\left( \frac{1}{%
d^{2}}\right) .$

What is more, in the complex-valued case there is an explicit formula for
all Lyapunov exponents even if $\Sigma $ is general. Namely, it is shown in 
\cite{forrester13} that 
\begin{equation}
\mu _{k}=\frac{1}{2}\Psi \left( k\right) +\frac{1}{2\prod\limits_{i<j}\left(
y_{i}-y_{j}\right) }\det \left[ 
\begin{array}{c}
\left[ y_{j}^{i-1}\right] _{i=1,\ldots ,k-1;j=1,\ldots ,d} \\ 
\left[ \left( \log y_{j}\right) y_{j}^{k-1}\right] _{j=1,\ldots ,d} \\ 
\left[ y_{j}^{i-1}\right] _{i=k+1,\ldots ,d;j=1,\ldots ,d}%
\end{array}%
\right] ,  \label{Forrester_main}
\end{equation}%
where $y_{j}$ are eigenvalues of $\Sigma ^{-1}.$ In particular for $k=1,$
one can re-write this as 
\begin{equation}
\mu _{1}=\frac{1}{2}\left[ \Psi \left( 1\right) -\sum_{j=1}^{d}\frac{\log
y_{j}}{\prod_{l\neq j}\left( 1-y_{j}/y_{l}\right) }\right]
\label{formula_largest_Forrester}
\end{equation}%
provided that all $y_{j}$ are different.

The proof of formula (\ref{Forrester_main}) is based on the
Harish-Chandra-Itzykson-Zuber integral and cannot be directly generalized to
the case of real or quaternion Gaussian matrices.

In fact, it appears that for the real-valued case with general $\Sigma ,$ an
explicit formula (due to Mannion \cite{mannion93}) is only known for
products of 2-by-2 Gaussian matrices:%
\begin{equation}
\mu _{1}=\frac{1}{2}\left[ \Psi \left( 1\right) +\log \left( \frac{1}{2}%
\mathrm{Tr}\Sigma +\sqrt{\det \Sigma }\right) \right] .  \label{Manning}
\end{equation}
(Some explicit formulas are also known for 2-by-2 random matrices with
non-Gaussian entries, see \cite{mtw08}. In addition, there are methods which
sometime allow one to compute Lyapunov exponents efficiently even when
explicit formulas are not available, see \cite{pollicott10a}.)

Our formula (\ref{formula_main3}) in Theorem \ref{theorem_main_formula}
provides an explicit formula applicable for real, complex and quaternion
matrices with general $\Sigma $. It is limited, however, to the case of the
largest Lyapunov exponent. For $\beta =2,$ Forrester's formula (\ref%
{formula_largest_Forrester}) can be derived from our formula by evaluating
the integral in (\ref{formula_main3}) using residues. Similarly, for $\beta
=4,$ one can use residues and derive the following formula:%
\begin{eqnarray}
2\mu _{1} &=&\Psi \left( 1\right) -\log \left( 2\right)
\label{quaternion_largest_Lyapunov} \\
&&-\left( \prod\limits_{i=1}^{d}y_{i}\right) ^{2}\sum_{i=1}^{d}\left\{ \frac{%
1}{y_{i}^{2}\prod\limits_{j\neq i}\left( y_{i}-y_{j}\right) ^{2}}\left[
1-\log y_{i}\left( 1+\sum_{j\neq i}\frac{2y_{i}}{y_{i}-y_{j}}\right) \right]
\right\} .  \notag
\end{eqnarray}%
In the remaining case of $\beta =1,$ the integral can be easily evaluated
numerically. In the asymptotic analysis of large-dimensional situations, we
will use our formula (\ref{formula_main3}) as most convenient even for $%
\beta =2$ and $4.$

\section{Proof of the formula for the largest Lyapunov exponent}

\label{section_integral_formula}

We start the proof of Theorem \ref{theorem_main_formula} by interpreting the
basic formula (\ref{basic_formula0}) for $k=1.$ Namely, by (\ref%
{basic_formula0}), $\mu _{1}$ is the expected logarithm of $\left\Vert
\Sigma ^{1/2}G_{1}\right\Vert $. Since $\left\Vert \Sigma
^{1/2}G_{1}\right\Vert ^{2}$ is a weighted sum of the squares of independent
Gaussian variables, its distribution is easy to calculate. Indeed, in the
real case $\beta =1,$ the characteristic function of the sum is 
\begin{equation*}
f\left( t\right) =\prod_{j=1}^{d}\left( 1-2i\sigma _{j}^{2}t\right) ^{-1/2}.
\end{equation*}%
(In this formula, the function $z^{-1/2}$ is determined by making a cut
along $z<0,$ and the selected branch of $\left( 1-2i\sigma _{j}^{2}t\right)
^{-1/2}$ equals $1$ at $t=0.$) Hence, the distribution density of the sum is 
\begin{eqnarray*}
p\left( \lambda \right) &=&\frac{1}{2\pi }\int_{-\infty }^{\infty
}e^{-i\lambda t}\prod_{j=1}^{d}\left( 1-2i\sigma _{j}^{2}t\right) ^{-1/2}dt
\\
&=&\frac{1}{2\pi }\left( \Pi _{j=1}^{d}y_{j}\right) ^{1/2}\int_{-\infty
}^{\infty }e^{-i\lambda t}\prod_{j=1}^{d}\left( y_{j}-2it\right) ^{-1/2}dt \\
&=&c\int_{0\mathcal{+}i\infty }^{0-i\infty }e^{-\lambda
z/2}\prod_{i=1}^{d}\left( z-y_{i}\right) ^{-1/2}dz,
\end{eqnarray*}%
where $y_{i}=1/\sigma _{i}^{2}$ and the normalization constant $c$ may
depend on $y_{i}$ and the choice of branch for functions $\left(
z-y_{i}\right) ^{-1/2}.$ By changing the contour, we get the expression 
\begin{equation*}
p\left( \lambda \right) =\frac{1}{c}\int_{\mathcal{C}}e^{-\lambda
z/2}\prod_{i=1}^{d}\left( z-y_{i}\right) ^{-1/2}dz,
\end{equation*}%
and the contour of integration goes in the counterclockwise direction around
points $y_{i}.$ In more detail, we can take the contour $\mathcal{C}$ that
starts at a large $R>\max \left\{ y_{i}\right\} ,$ goes along the upper edge
of the real axis to $r>0,$ which is smaller than all of $y_{i},$ then
crosses to the lower edge of the real axis, and then returns along this
lower edge to $R.$ We then take the limit for $R\rightarrow \infty .$

For general $\beta $ (i.e., $\beta =1,2,$ or $4$), we can similarly obtain 
\begin{equation*}
p_{\beta }\left( \lambda \right) =\frac{1}{c_{\beta }}\int_{\mathcal{C}%
}e^{-\beta \lambda z/2}\prod_{i=1}^{d}\left( z-y_{i}\right) ^{-\beta /2}dz.
\end{equation*}

Next, we use formula (\ref{basic_formula0}) and find that 
\begin{equation}
\mu _{1}=\frac{1}{2c_{\beta }}\int_{0}^{\infty }\log \lambda \left[ \int_{%
\mathcal{C}}e^{-\beta \lambda z/2}\prod_{i=1}^{d}\left( z-y_{i}\right)
^{-\beta /2}dz\right] d\lambda ,  \label{formula_main}
\end{equation}%
where%
\begin{equation}
c_{\beta }=\int_{0}^{\infty }\int_{\mathcal{C}}e^{-\beta \lambda
z/2}\prod_{i=1}^{d}\left( z-y_{i}\right) ^{-\beta /2}dzd\lambda .
\label{formula_main_constant}
\end{equation}

By changing the order of integration in (\ref{formula_main}) and computing
the inner integral, we find (with the help of the identity $\int_{0}^{\infty
}\left( \log t\right) e^{-t}dt=-\gamma $) that 
\begin{equation}
2\mu _{1}=-\gamma +\frac{1}{\widetilde{c}_{\beta }}\left[ \frac{1}{2\pi i}%
\int_{\mathcal{C}}\left( \frac{2}{\beta z}\right) \log \left( \frac{2}{\beta
z}\right) \prod_{i=1}^{d}\left( z-y_{i}\right) ^{-\beta /2}dz\right] ,
\label{formula_main2a}
\end{equation}%
where 
\begin{eqnarray*}
\widetilde{c}_{\beta } &=&\frac{c_{\beta }}{2\pi i}=\frac{1}{2\pi i}\int_{%
\mathcal{C}}\left( \frac{2}{\beta z}\right) \prod_{i=1}^{d}\left(
z-y_{i}\right) ^{-\beta /2}dz \\
&=&\frac{2}{\beta }\prod\limits_{i=1}^{d}\left( -y_{i}\right) ^{-\beta /2}.
\end{eqnarray*}%
This implies that%
\begin{equation}
2\mu _{1}=\Psi \left( 1\right) +\log \left( \frac{2}{\beta }\right) +\frac{1%
}{2\pi i}\int_{\mathcal{C}}\log \left( z\right) \prod_{i=1}^{d}\left( 1-%
\frac{z}{y_{i}}\right) ^{-\beta /2}\frac{dz}{z}.  \label{formula_main2}
\end{equation}

\begin{figure}[tbph]
\includegraphics[width=12cm]{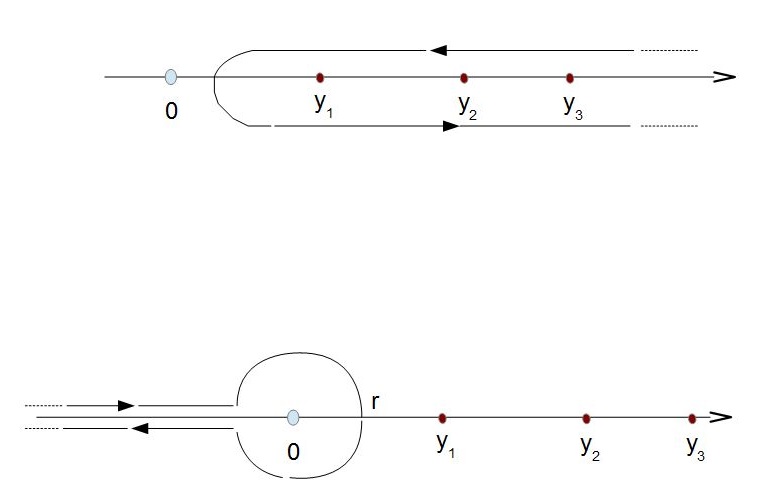} \label{fig:contour}
\caption{Original and modified contours.}
\end{figure}

Now, let us change the contour. As before, let the contour start at a large $%
R>\max \left\{ y_{i}\right\} ,$ goes along the upper edge of the real axis
to $r>0,$ which is smaller than all of $y_{i},$ then crosses to the lower
edge of the real axis, and then returns along this lower edge to $R.$ 

Next, let us move this contour of integration so that the new contour $%
\mathcal{C}^{\prime }$ starts at $-R,$ goes along the upper edge of the real
axis to $-r,$ then goes along the circle of radius $r$ around the $0$ in the
clockwise direction, and then returns to $-R$ along the lower edge of the
real axis. The change in the integral corresponding to the change from
contour $\mathcal{C}$ to contour $\mathcal{C}^{\prime }$ is small for large $%
R.$

Indeed, consider the contour that starts at $-R,$ then goes to $-r$ along
the upper edge of the real axis, then to $r$ clockwise along the circle
centered at zero, and then to $R$ along the upper edge of the real axis.
Finally its returns back to $-R$ conterclocwise along the circle centered at
zero. The integral over this closed contour is zero because the integrand is
holomorphic inside the contour. In addition, the integral over the
semicircle with radius $R$ is small for large $R.$ Hence, the integral over
the first portion of this contour is also small. This shows that the
difference in integrals is small for those portions of the contours $%
\mathcal{C}$ and $\mathcal{C}^{\prime }$ that are above the real axis.
Similar we can treat the portions of the contours $\mathcal{C}$ and $%
\mathcal{C}^{\prime }$ that are below the real axis.

Then by computing the integrals over the two rays and the circle, and by
letting  $r\rightarrow 0,$ we find that 
\begin{eqnarray*}
\lim_{\substack{ R\rightarrow \infty , \\ r\rightarrow 0}}\frac{-1}{2\pi i}%
\int_{\mathcal{C}^{\prime }}\log \left( z\right) \prod_{i=1}^{d}\left( 1-%
\frac{z}{y_{i}}\right) ^{-\beta /2}\frac{dz}{z} &=&\lim_{r\rightarrow
0}\left\{ \int_{r}^{\infty }\prod_{i=1}^{d}\left( 1+\frac{x}{y_{i}}\right)
^{-\beta /2}\frac{dx}{x}+\log r\right\}  \\
&=&\int_{0}^{1}\left( \prod_{i=1}^{d}\left( 1+\frac{x}{y_{i}}\right)
^{-\beta /2}-1\right) \frac{dx}{x} \\
&&+\int_{1}^{\infty }\prod_{i=1}^{d}\left( 1+\frac{x}{y_{i}}\right) ^{-\beta
/2}\frac{dx}{x}.
\end{eqnarray*}%
(In the first integral we used the fact that the branches of the logarithm
on the upper and the lower edges of the real axis differ by $2\pi i$.) This
proves formula (\ref{formula_main3}). \hfill $\square $

\section{Asymptotic behavior for large matrix size}

\label{section_spike}

\subsection{A model without spikes}

We start the proof of Theorem \ref{theorem_no_spikes} by noting that a
straightforward calculation gives the identity

\begin{equation*}
\Psi \left( 1\right) +\log \left( \frac{2}{\beta }\right) +\int_{0}^{\infty }%
\left[ 1_{\left[ 0,1\right] }\left( x\right) -e^{-\left( \beta /2\right)
\lambda x}\right] \frac{dx}{x}=\log \lambda .
\end{equation*}%
If we compare this with the formula (\ref{formula_main3}), then we find that
it is sufficient to show that 
\begin{equation}
\int_{0}^{\infty }\left[ \prod_{i=1}^{d}\left( 1+\frac{\theta _{i}}{d}%
x\right) ^{-\beta /2}-e^{-\left( \beta /2\right) \lambda x}\right] \frac{dx}{%
x}\rightarrow 0  \label{convergence}
\end{equation}%
as $d$ $\rightarrow \infty .$ First, consider 
\begin{equation*}
I_{1}:=\int_{0}^{\infty }\left[ \prod_{i=1}^{d}\left( 1+\frac{\theta _{i}}{d}%
x\right) ^{-\beta /2}-\prod_{i=1}^{d}\exp \left( -\frac{\beta }{2}\frac{%
\theta _{i}}{d}x\right) \right] \frac{dx}{x}.
\end{equation*}

We split $I_{1}$ in two integrals, $I_{1}^{\prime }+I_{1}^{\prime \prime },$
the first one is over the interval from \thinspace $0$ to $M,$ and the other
is over the interval from $M\,\ $to infinity. For the second integral, we
have 
\begin{equation*}
\int_{M}^{\infty }\exp \left( -\frac{\beta }{2}\frac{\sum \theta _{i}}{d}%
x\right) \frac{dx}{x}=\int_{M\left( \frac{\beta }{2}\right) \frac{1}{d}\sum
\theta _{i}}^{\infty }\frac{e^{-t}dt}{t}\leq e^{-M\left( \beta /2\right)
}\rightarrow 0
\end{equation*}%
as $M\rightarrow \infty $ and the convergence is uniform in $d.$ In
addition, 
\begin{eqnarray*}
\int_{M}^{\infty }\prod_{i=1}^{d}\left( 1+\frac{\theta _{i}}{d}x\right)
^{-\beta /2}\frac{dx}{x} &\leq &\int_{M/d}^{\infty }\left( 1+t\right)
^{-d\beta /2}\frac{dt}{t} \\
&\leq &\frac{d}{M}\int_{M/d}^{\infty }\left( 1+t\right) ^{-d\beta /2}dt \\
&=&\frac{d}{M}\frac{1}{d\left( \beta /2\right) -1}\left( 1+\frac{M}{d}%
\right) ^{-d\beta /2+1}\rightarrow 0
\end{eqnarray*}%
as $M\rightarrow \infty ,$ again uniformly in $d>1.$ We conclude that for
every $\varepsilon >0,$ we can find $M_{0}$ such that $\left\vert
I_{1}^{\prime \prime }\right\vert \leq \varepsilon $ for all $M\geq M_{0}$
and all $d>1.$ In words, we can make $I_{1}^{\prime \prime }$ arbitrarily
small uniformly in $d$ by taking $M$ sufficiently large.

For the integral $I_{1}^{\prime }$, we estimate the integrand by using the
fact that if $\left\vert z_{i}\right\vert \leq 1$ and $\left\vert
w_{i}\right\vert \leq 1,$ then 
\begin{equation*}
\left\vert \prod_{i=1}^{d}z_{i}-\prod_{i=1}^{d}w_{i}\right\vert \leq
\sum_{i=1}^{d}\left\vert z_{i}-w_{i}\right\vert ,
\end{equation*}%
(see Lemma 1 of Section 27 in Billingsley \cite{billingsley95}). Since 
\begin{equation*}
\left\vert \left( 1+\frac{\theta _{i}}{d}x\right) ^{-\beta /2}-\exp \left( -%
\frac{\beta }{2}\frac{\theta _{i}}{d}x\right) \right\vert \leq C\left( \frac{%
\theta _{i}}{d}x\right) ^{2}
\end{equation*}%
for all $x\leq d/L,$ therefore (for $d\geq LM$), we estimate 
\begin{equation*}
I_{1}^{\prime }\leq \int_{0}^{M}C\frac{L^{2}}{d}xdx=\frac{CL^{2}M^{2}}{2d}.
\end{equation*}%
For a fixed $M,$ this can be made arbitrarily small by choosing $d$
sufficiently large. Hence, $I_{1}\rightarrow 0$ as $d\rightarrow \infty .$

Similarly, 
\begin{equation*}
I_{2}^{\prime \prime }:=\int_{M}^{\infty }\left( e^{-\left( \beta /2\right)
\lambda x}-e^{-\left( \beta /2\right) d^{-1}\left( \sum \theta _{i}\right)
x}\right) \frac{dx}{x}\rightarrow 0
\end{equation*}%
as $M\rightarrow \infty $ uniformly in $d,$ and for 
\begin{equation*}
I_{2}^{\prime }:=\int_{0}^{M}\left( e^{-\left( \beta /2\right) \lambda
x}-e^{-\left( \beta /2\right) d^{-1}\left( \sum \theta _{i}\right) x}\right) 
\frac{dx}{x},
\end{equation*}%
we estimate 
\begin{eqnarray*}
\left\vert I_{2}^{\prime }\right\vert &\leq &\int_{0}^{M}C\left\vert \lambda
-d^{-1}\left( \sum \theta _{i}\right) \right\vert dx \\
&=&CM\left\vert \lambda -d^{-1}\left( \sum \theta _{i}\right) \right\vert
\rightarrow 0
\end{eqnarray*}%
as $d\rightarrow \infty $ for a fixed $M.$ Altogether, the convergence of $%
I_{1}$ and $I_{2}:=I_{2}^{\prime }+I_{2}^{^{\prime \prime }}$ to zero proves
(\ref{convergence}) and completes the proof. \hfill $\square $

\subsection{A model with a spike}

In order to prove Theorem \ref{theorem_spikes}, consider a slightly modified
model in which the eigenvalues of $\Sigma _{d}$ are $\lambda _{i}^{\left(
d\right) }=1$ for $i=1,\ldots ,d-1,$ and $\lambda _{d}^{\left( d\right)
}=\theta =d/t>1$. This model is obtained from the model in Theorem \ref%
{theorem_spikes} by multiplying $\Sigma _{d}$ by $d.$ Therefore, the results
in the theorem can be recovered by subtracting $\frac{1}{2}\log d$ from the
largest Lyapunov exponent of the modified model.

For the complex case, we have 
\begin{eqnarray*}
2\mu _{1} &=&\Psi \left( d\right) +\int_{0}^{\infty }\frac{1}{\left(
1+x\right) ^{d-1}}\left[ \frac{1}{1+x}-\frac{1}{1+\theta x}\right] \frac{dx}{%
x} \\
&=&\Psi \left( d\right) +\int_{0}^{\infty }\frac{1}{\left( 1+u/d\right) ^{d}}%
\left[ \frac{1}{1+u/d}-\frac{1}{1+u/t}\right] \frac{du}{u} \\
&=&\log d+\int_{0}^{\infty }e^{-u}\left[ 1-\frac{1}{1+u/t}\right] \frac{du}{u%
}+O_{t}\left( 1/d\right) \\
&=&\log d+\int_{0}^{\infty }e^{-xt}\frac{dx}{1+x}+O_{t}\left( 1/d\right) \\
&=&\log d+e^{t}\int_{1}^{\infty }e^{-xt}\frac{dx}{x}+O_{t}\left( 1/d\right) .
\end{eqnarray*}

For the real case, we assume $d=2k$ (the other case is similar) and write:%
\begin{eqnarray*}
2\mu _{1} &=&\log 2+\Psi \left( k\right) +\int_{0}^{\infty }\frac{1}{\left(
1+x\right) ^{k-1/2}}\left[ \frac{1}{\sqrt{1+x}}-\frac{1}{\sqrt{1+\theta x}}%
\right] \frac{dx}{x} \\
&=&\log 2+\Psi \left( k\right) +\int_{0}^{\infty }\frac{1}{\left( 1+\frac{u}{%
k}\right) ^{k-1/2}}\left[ \frac{1}{\sqrt{1+u/k}}-\frac{1}{\sqrt{1+2u/t}}%
\right] \frac{du}{u} \\
&=&\log d+\int_{0}^{\infty }e^{-u}\left[ 1-\frac{1}{\sqrt{1+2u/t}}\right] 
\frac{du}{u}+O_{t}\left( 1/d\right) \\
&=&\log d+e^{t/2}\int_{1}^{\infty }e^{-tx/2}\frac{dx}{\sqrt{x}\left( \sqrt{x}%
+1\right) }+O_{t}\left( 1/d\right) .
\end{eqnarray*}%
(The last step uses the change of variable $x=1+2u/t.$) \hfill $\square $

\textbf{Proof of Theorem \ref{theorem_no_spikes_2nd}:} Again, it is more
convenient to use the model in which the eigenvalues of $\Sigma _{d}$ are $%
\lambda _{i}^{\left( d\right) }=1$ for $i=1,\ldots ,d-1,$ and $\lambda
_{d}^{\left( d\right) }=\theta >1$. It is elementary to compute that 
\begin{equation*}
\int_{0}^{\infty }\left[ \boldsymbol{1}_{[0,1]}\left( x\right) -\left(
1+x\right) ^{-d}\right] \frac{dx}{x}=1+\frac{1}{2}+\ldots +\frac{1}{d-1}%
=\Psi \left( d\right) -\Psi (1).
\end{equation*}%
(For example, one can use the identity $\frac{1}{\left( 1+x\right) ^{d}x}=%
\frac{1}{x}-\frac{1}{1+x}-\ldots -\frac{1}{\left( 1+x\right) ^{d}}.$) Hence,
by using formula (\ref{formula_main3}), we can write 
\begin{align}
2\mu _{1}& =\Psi \left( 1\right) +\int_{0}^{\infty }\left[ \boldsymbol{1}%
_{[0,1]}\left( x\right) -\frac{1}{\left( 1+x\right) ^{d-1}\left( 1+\theta
x\right) }\right] \frac{dx}{x}  \notag \\
& =\Psi \left( d\right) +\int_{0}^{\infty }\left[ \frac{1}{\left( 1+x\right)
^{d-1}}\left( \frac{1}{1+x}-\frac{1}{1+\theta x}\right) \right] \frac{dx}{x},
\notag \\
& =\Psi \left( d\right) +f_{d},  \label{formula_mu_fd}
\end{align}%
where 
\begin{equation}
f_{d}:=\left( \theta -1\right) \int_{0}^{\infty }\frac{1}{\left( 1+x\right)
^{d}\left( 1+\theta x\right) }dx\leq \frac{\theta -1}{d}
\label{def_additional_term}
\end{equation}%
One can check that for $d\geq 1,$ the additional term $f_{d}$ satisfies the
recursion equation:%
\begin{equation}
f_{d}=\left( \frac{\theta -1}{\theta }\right) \left( \frac{1}{d}%
+f_{d+1}\right) .  \label{recursion}
\end{equation}%
Hence, we obtain a convergent series for $f_{d},$%
\begin{equation*}
f_{d}=s\sum_{k=0}^{\infty }\frac{s^{k}}{d+k},
\end{equation*}%
where $s:=\left( \theta -1\right) /\theta <1.$ Hence, 
\begin{eqnarray*}
\left\vert df_{d}-\frac{s}{s-1}\right\vert &=&s\left\vert \sum_{k=0}^{\infty
}\left( \frac{s^{k}}{1+k/d}-s^{k}\right) \right\vert \\
&=&\frac{s}{d}\sum_{k=0}^{\infty }\frac{ks^{k}}{1+k/d}\leq \frac{1}{d}\left( 
\frac{s}{s-1}\right) ^{2}.
\end{eqnarray*}%
That is, 
\begin{equation*}
\left\vert f_{d}-\frac{\theta -1}{d}\right\vert \leq \left( \frac{\theta -1}{%
d}\right) ^{2}.
\end{equation*}%
Together with the asymptotic expansion for the digamma function, $\Psi
\left( d\right) =\log d-\frac{1}{2d}+O\left( \frac{1}{d^{2}}\right) ,$ this
limit implies the statement of the theorem. \hfill $\square $

\textbf{Proof of Theorem \ref{theorem_spikes_2nd}:} We can use the recursion
in (\ref{recursion}) and the initial condition $f_{1}=\log \theta $ in order
to obtain 
\begin{equation}
f_{d}=\left( \frac{\theta }{\theta -1}\right) ^{d-1}\left( \log \theta
-\sum_{k=1}^{d-1}\frac{1}{k}\left( 1-\frac{1}{\theta }\right) ^{k}\right) .
\label{fd_formula2}
\end{equation}

Note that 
\begin{eqnarray*}
\left\vert \left( \frac{\theta }{\theta -1}\right) ^{d-1}-1\right\vert
&=&\left\vert \exp \left[ \left( d-1\right) \log \left( 1+\frac{1}{\theta -1}%
\right) \right] -1\right\vert \\
&=&\left\vert \exp \left[ \frac{d-1}{\theta -1}+O\left( \frac{d-1}{\left(
\theta -1\right) ^{2}}\right) \right] -1\right\vert \\
&=&O\left( \frac{d-1}{\theta -1}\right) .
\end{eqnarray*}%
Similarly, we estimate 
\begin{equation*}
\sum_{k=1}^{d-1}\frac{1}{k}\left( 1-\frac{1}{\theta }\right)
^{k}=\sum_{k=1}^{d-1}\frac{1}{k}\left[ 1+O\left( \frac{d}{\theta }\right) %
\right] .
\end{equation*}%
Since $\sum_{k=1}^{d-1}\frac{1}{k}=\Psi \left( d\right) -\Psi \left(
1\right) \sim \log d,$ we find that 
\begin{equation*}
f_{d}=\log \theta -\Psi \left( d\right) +\Psi \left( 1\right) +O\left( \frac{%
d}{\theta }\log \theta \right) .
\end{equation*}

Therefore by (\ref{formula_mu_fd}), $2\mu _{1}=\log \theta -\gamma +O\left( 
\frac{d}{\theta }\log \theta \right) .$ \hfill $\square $

\section{Conclusion}

\label{section_conclusion}

We derived a new formula for the largest Lyapunov exponent of Gaussian
matrices and the asymptotic expressions for this exponent when the matrices
are large. The asymptotic expressions are derived for two distinct cases,
one is when the covariance matrix has no spikes and another one is when it
has a single spike.

\begin{figure}[tbph]
\includegraphics[width=12cm]{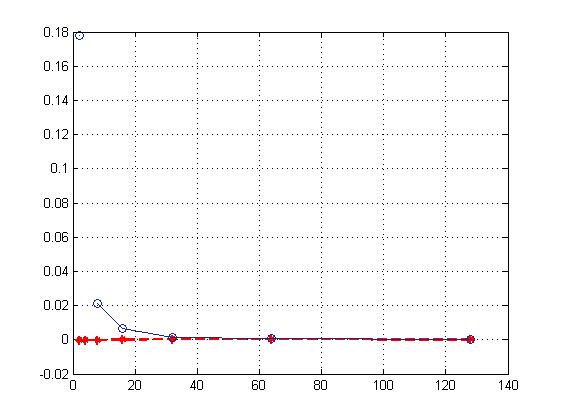} \label{fig:GB}
\caption{The deviation of the largest Lyapunov exponent from the theoretical prediction. The no-spike case. The dashed line is for Gaussian and the solid line is for 
Bernoulli matrices. The horizontal axis shows the size of
the matrix.}
\end{figure}
\begin{figure}[tbph]
\includegraphics[width=12cm]{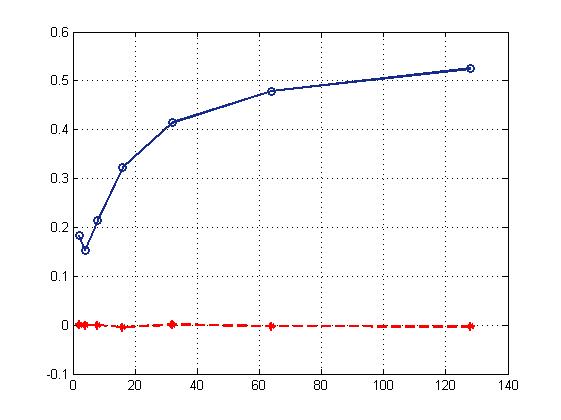} \label{fig:GB_spike}
\caption{The deviation of the largest Lyapunov exponent from the theoretical prediction. The spike case.}
\end{figure}

An interesting question is whether the asymptotic expressions derived in
this paper remain valid for non-Gaussian random matrices with independent
entries. Some numerical evidence concerning this question is presented in
Figures 3 and 4.

In order to create plots in these figures, we estimate the Lyapunov exponent
using simulations and then plot its difference from the theoretical
prediction in the formula (\ref{formula_main3}). The difference is shown
against the dimension of the matrix in order to check whether it declines
with the increase in dimension as the concept of universality would predict.

Figure 3 presents this plot for the matrices $GA,$ where $G$ is
the diagonal matrix with half of the diagonal entries equal to $1$ and the
remaining half equal to $2.$ The dashed line is for a $d$-by-$d$ Gaussian
matrix $A$ where each entry is independent and have variance $d^{-1}.$ It is
very close to zero even for relatively small $d.$

The solid line is for a $d$-by-$d$ matrix $A$ where each entry is
independent and equals $\pm d^{-1/2}$ with probability $1/2.$ The solid line
approaches zero as $d$ grows, which suggests that the theoretical
approximation gives a good approximation for the Lyapunov exponent of the
Bernoulli matrix $A.$

(This numerical evidence should be taken with a grain of salt. While it
might indicate the universality of the theoretical approximation for a large
class of matrices, in fact this approximation is not valid for Bernoulli
matrices (at least in the strict sense). Indeed, with positive probability
the Bernoulli matrix $A$ will have all entries equal to $d^{-1/2}$ and
therefore the starting vector will be mapped to a vector proportional to $%
\left( 1,\ldots ,1,2,\ldots ,2\right) .$ For $d\geq 4,$ this vector has a
positive probability to be mapped to $0$ by an appropriate Bernoulli matrix $%
A.$ Hence, with probability $1$ the product of independent matrices $A_{k}$
will eventually maps an arbitrary starting vector to zero. It can be seen
that this happens for $d=4$ in our plot where the data point is missing.
However, for large $d$ the simulations would need a prohibitively large time
to detect this phenomenon.)

Figure 4 shows an analogous plot for the case when the
matrix $G$ has a spike. Namely, $G$ is assumed to be diagonal with all
diagonal entries except one equal to $1$ and the exceptional diagonal entry
equal to $d.$ This plot shows that the theoretical approximation from
Theorem \ref{theorem_main_formula} is valid for Gaussian but invalid for
Bernoulli matrices. This raises a natural question of how the approximation
should be modified for large non-Gaussian matrices with a spike in the
covariance matrix.

\bigskip

\bibliographystyle{plain}
\bibliography{comtest}

\begin{thebibliography}{10}

\bibitem{billingsley95}
Patrick Billingsley.
\newblock {\em Probability and Measure}.
\newblock John Wiley and Sons, third edition, 1995.

\bibitem{bougerol_lacroix85}
P.~Bougerol and J.~Lacroix.
\newblock {\em Products of random matrices with applications to Schrodinger
  operators}, volume~8 of {\em Progress in probability and statistics}.
\newblock Birkhauser: Boston, 1985.

\bibitem{cohen_kesten_newman86}
Joel~E. Cohen, Harry Kesten, and Charles~M. Newman, editors.
\newblock {\em Random Matrices and Their Applications}, volume~50 of {\em
  Contemporary Mathematics}.
\newblock American Mathematical Society, 1986.

\bibitem{cohen_newman84}
Joel~E. Cohen and Charles~M. Newman.
\newblock The stability of large random matrices and their products.
\newblock {\em Annals of Probability}, 12:283--310, 1984.

\bibitem{forrester13}
P.~J. Forrester.
\newblock Lyapunov exponents for products of complex \mbox{G}aussian random
  matrices.
\newblock {\em Journal of Statistical Physics}, 151:796--808, 2013.
\newblock \mbox{\href{http://arxiv.org/abs/1206.2001}{arxiv:1206.2001}}.

\bibitem{furstenberg_kesten60}
H.~Furstenberg and H.~Kesten.
\newblock Products of random matrices.
\newblock {\em Annals of Mathematical Statistics}, 31:457--469, 1960.

\bibitem{furstenberg63}
Harry Furstenberg.
\newblock Noncommuting random products.
\newblock {\em Transactions of the American Mathematical Society},
  108:377--428, 1963.

\bibitem{guivarch_raugi85}
Y.~Guivarch and A.~Raugi.
\newblock Frontiere de \mbox{F}urstenberg, proprietes de contraction et
  theoremes de convergence.
\newblock {\em Zeit. Fur Wahrscheinlichkeitstheorie und Verw. Gebiete},
  67:265--278, 1985.

\bibitem{kargin07}
V.~Kargin.
\newblock The norm of products of free random variables.
\newblock {\em Probability Theory and Related Fields}, 139:397--413, 2007.
\newblock \mbox{\href{http://arxiv.org/abs/math/0611593}{arxiv:math/0611593}}.

\bibitem{kargin08}
V.~Kargin.
\newblock Lyapunov exponents of free operators.
\newblock {\em Journal of Functional Analysis}, 255:1874--1888, 2008.
\newblock \mbox{\href{http://arxiv.org/abs/0712.1378}{arxiv:0712.1378}}.

\bibitem{lepage89}
Emile le~Page.
\newblock R\mbox{\'{e}}gularit\mbox{\'{e}} de plus grand exposant
  charact\mbox{\'{e}}ristique de produit de matrices al\mbox{\'{e}}atoires
  ind\mbox{\'{e}}pendantes et applications.
\newblock {\em Annales de l'Institut Henri Poincaré. Probabilités et
  Statistiques}, 25:109--142, 1989.

\bibitem{mannion93}
David Mannion.
\newblock Products of $2x2$ random matrices.
\newblock {\em Annals of Applied Probability}, 3:1189--1218, 1993.

\bibitem{mtw08}
Jens Marklof, Yves Tourigny, and Lech Wolowski.
\newblock Explicit invariant measures for products of random matrices.
\newblock {\em Transactions of the American Mathematical Society},
  360:3391--3427, 2008.

\bibitem{newman86a}
C.~M. Newman.
\newblock Lyapunov exponents for some products of random matrices: Exact
  expressions and asymptotic distributions.
\newblock In Joel~E. Cohen, Harry Kesten, and Charles~M. Newman, editors, {\em
  Random Matrices and Their Applications}, volume~50 of {\em Contemporary
  Mathematics}, pages 183--195. American Mathematical Society, 1986.

\bibitem{newman86b}
Charles~M. Newman.
\newblock The distribution of \mbox{L}yapunov exponents: Exact results for
  random matrices.
\newblock {\em Communications in Mathematical Physics}, 103:121--126, 1986.

\bibitem{oseledec68}
V.~I. Oseledec.
\newblock A multiplicative ergodic theorem. \mbox{L}japunov characteristic
  numbers for dynamical systems.
\newblock {\em Transactions of the Moscow Mathematical Society}, 19:197--231,
  1968.

\bibitem{peres92}
Yuval Peres.
\newblock Domains of analytic continuation for the top \mbox{L}yapunov
  exponent.
\newblock {\em Ann. I. H. Poincare}, 28:131--148, 1992.

\bibitem{pollicott10a}
Mark Pollicott.
\newblock Maximal \mbox{L}yapunov exponents for random matrix products.
\newblock {\em Inventiones Mathematicae}, 181:209--226, 2010.

\bibitem{raghunathan79}
M.~S. Raghunathan.
\newblock A proof of \mbox{O}seledec\mbox{'}s multiplicative ergodic theorem.
\newblock {\em Israel Journal of Mathematics}, 32:356--362, 1979.

\bibitem{ruelle79}
David Ruelle.
\newblock Analyticity properties of the characterisic exponents of random
  matrix products.
\newblock {\em The Advances of Mathematics}, 32:68--80, 1979.

\bibitem{ruelle82}
David Ruelle.
\newblock Characterisic exponents and invariant manifolds in \mbox{H}ilbert
  space.
\newblock {\em The Annals of Mathematics}, 115:243--290, 1982.

\bibitem{ruelle84}
David Ruelle.
\newblock Characterisic exponents for a viscous fluid subjected to time
  dependent forces.
\newblock {\em Communications in Mathematical Physics}, 93:285--300, 1984.

\end{thebibliography}

\end{document}